\documentclass[reqno,b5paper]{amsart}

\usepackage{amsmath}
\usepackage{amssymb}
\usepackage{amsthm}
\usepackage{enumerate}
\usepackage[mathscr]{eucal}
\usepackage{eqlist}

\theoremstyle{plain}
\newtheorem{theorem}{Theorem}[section]

\theoremstyle{definition}
\newtheorem{definition}{Definition}[section]
\theoremstyle{proposition}
\newtheorem{proposition}{Proposition}[section]
\newtheorem{remark}{\textnormal{\textbf{Remark}}}

\theoremstyle{remark}

\numberwithin{equation}{section}

\begin{document}
\title[Diophantine Equation $X^4+Y^4=2(U^4+V^4)$]%
{Diophantine Equation $X^4+Y^4=2(U^4+V^4)$}
\author[Farzali Izadi \and Kamran Nabardi]%
{Farzali Izadi* \and Kamran Nabardi**}

\newcommand{\acr}{\newline\indent}

\address{\llap{*\,} Department of Mathematics \acr
                   Azarbaijan Shahid Madani University\acr
                   Tabriz, Iran}
\email{izadi@azaruniv.edu}

\address{\llap{**\,}Department of Mathematics\acr
                    Azarbaijan Shahid Madani University\acr
                    Tabriz, Iran}
\email{nabardi@azaruniv.edu}



\begin{abstract}
In this paper, the theory of elliptic curves is used for finding the solutions of the quartic Diophantine equation $X^4+Y^4=2(U^4+V^4)$.
\end{abstract}

\maketitle

\section{Itroduction}
The solubility in integers  of diophantine equation
\begin{equation}\label{EQ1}
a_0{X}^4+a_1{Y}^4+a_2{U}^4+a_3{V}^4=0
\end{equation}
 has been considered by many mathematicians, where $a_0, a_1, a_2$ and $a_3$ are nonzero integers. The most famous and simplest one is proposed by Euler  (see \cite {Ha} page 201) for the constants $a_0=a_1=1$ and $a_2=a_3=-1$. Euler gave a two-parameter solutions for this equation.
 Zajta  \cite{Zaj} applied several methods including the Pythagorean and algebraic reduction method, for parametrization of $A^4+B^4=C^4+D^4$.
  Brudno \cite{Bru} and Lander \cite{Lan} gave new parametrizations for like wise power diophantine equations, specially for $A^4+B^4=C^4+D^4$.
  Using geometric methods and the property of tangent plane, Richmond \cite{Rich} parameterized the equation \eqref{EQ1}, for the case that  the product $a_0a_1a_2a_3$ is a square number.
 Setting $a_0=1, a_1=4, a_2=-1$ and $a_3=-4$, Choudhry \cite{Chou} presented two-parameter solutions of the equation .
 Choudhry \cite{Chou2},  has considered a special family of Diophantine equation by
 \begin{equation}
 A^4+hB^4=C^4+hD^4,
 \end{equation}
 and has found a list of integer solutions for  cases $h\leq 101$.
  Noam Elkies  \cite{Elk} found infinitely many solutions of equation \eqref{EQ1}
 by taking $a_0=a_1=a_2=1$ and $a_3=-1$.  In his method he used the theory of elliptic curves .
This paper is concerned with the integral solutions of \eqref{EQ1} where $a_0=a_1=1$ and $a_2=a_3=-2$, i.e.,
\begin{equation}\label{EQ2}
X^4+Y^4=2(U^4+V^4).
\end{equation}
The smallest known solution for this equation is $(X, Y, U, V)=(21,19,20,7)$. When we say the smallest solution we mean the  smallest up to sign. For example $(21,19,20,-7)$ is a solution but it is not a new one. We give infinitely many solutions of \eqref{EQ2} by means of a specific congruent number elliptic curve namely  $y^2=x^3-36x^2$. \\
  First, let us recall some basic facts about elliptic curves. An elliptic curve $E$ over  $\Bbb{Q}$ is a curve  that is given by an equation of the form
$y^2=x^3+ax+b$, where $a, b\in \Bbb{Q}.$
By the Mordell-Weil theorem, the rational points on an elliptic curve form a finitely  generated abelian group, which is denoted by $E(\Bbb{Q})$ and so one can write the following decomposition
\begin{equation}
E(\Bbb{Q})\simeq E(\Bbb{Q})_{{\rm{tors}}}\oplus \Bbb{Z}^r,
\end{equation}
where $r$ is a nonnegative  integer called the  $rank$ of $E$ and $E(\Bbb{Q})_{{\rm{tors}}}$ is the finite group
 consisting of all the elements of finite order in $E(\Bbb{Q})$ \cite{Sil}.\\
A positive square free  integer   $n$ is called a congruent number if it is the  area of some right triangle with rational sides. The following theorem tells us that whether a number is congruent or not.
\begin{theorem}\label{th1}
Consider the elliptic curve $E_n: y^2=x^3-n^2x$. $n$ is a congruent number if and only if the elliptic curve $E_n(\Bbb{Q})$ has a positive rank.
\end{theorem}
\begin{proof}
See \cite{Kob}
\end{proof}
In  \eqref{EQ2}, let us change $X$ to $U+t$ and $Y$ to $U-t$ where $t$ is a parameter. Therefore we have
the equation $6t^2U^2+t^4=V^4$. Now  taking $Z=tU$ yields
\begin{equation}\label{EQ3}
6Z^2=V^4-t^4.
\end{equation}
Having said that, the following proposition is useful for our purpose.

\begin{proposition}\label{pro1}
Let $c$ be a nonzero integer. The equation $X^4-Y^4=cZ^2$ has a solution with $XYZ\not=0$ if and only if $|c|$ is a congruent number. More precisely, if $X^4-Y^4=cZ^2$ with $XYZ\not=0$ then $E_c: y^2=x^3-c^2x$ with $(x,y)=(-cY^2/X^2,c^2YZ/X^3)$, and conversely if $E_c: y^2=x^3-c^2x$ with $y\not=0$ then $X^4-Y^4=cZ^2$, with\\
\begin{equation*}
\begin{array}{llll}
X=x^2+cx-c^2,&Y=x^2-2cx-c^2,&\text{and}& Z=4y(x^2+c^2).
\end{array}
\end{equation*}
\end{proposition}
\begin{proof}
 See section 6.5 proposition 6.5.6 of \cite{Coh}.
\end{proof}
According to the equation  \eqref{EQ2},  theorem \ref{th1}  and  proposition \ref{pro1}, we see that the equation \eqref{EQ3} has a solution, since $c=6$ is a congruent number. So, equation \eqref{EQ2} has a solution. To find this solution we use the transformations of proposition \ref{pro1}. From $(x,y)$ on elliptic curve $E_6=x^3-36x$, we obtain
\begin{equation*}
\begin{array}{l}
t=x^2-12x-36,\\
\\
V=x^2+12x-36,\\
\\
U={\frac{4y(x^2+36)}{t}}.
\end{array}
\end{equation*}
Therefore, $(U+t, U-t, U, V)$ is a rational solution of \eqref{EQ2}. Multiplying this solution by $t$, we eliminate  denominator  of $U$. Next, we let $x=b/e^2$ and $y=c/e^3$ for some integers $b, c, e.$ Substituting these $x$ and $y$ and multiplying all equations by $e^8$ we get the following integer solutions for \eqref{EQ2}.
\begin{equation}\label {EQ4}
\begin{array}{l}
X=b^4 + 1296e^8 + 864be^6 + 72b^2e^4 + 144ce^5 - 24b^3e^2 + 4b^2ce,\\
\\
Y=-864be^6 - b^4 - 1296e^8 - 72b^2e^4 + 144ce^5 + 24b^3e^2 + 4b^2ce,\\
\\
U=4(b^2 + 36e^4)ce,\\
\\
V=(b^2 - 36e^4 - 12be^2)(b^2 - 36e^4 + 12be^2).\\
\end{array}
\end{equation}
\begin{remark}
 Note that, the additive inverse of a point $(x,y)$ on $E_6$ is $(x,-y)$. This means that we change $c$ to $-c$ in \eqref{EQ4}.  Consequently, If $(X,Y,U,V)$ is a
 solution obtained from $(x,y)$, then $(-X,-Y,-U,V)$  is a solution obtained from $(x,-y)$, which is not a new one up to sign.
 \end{remark}
\section{Numerical Results}
In this section we obtain primitive solutions of the diophantine equation \eqref{EQ2}.
\begin{definition}
 A solution $(A,B,C,D)$ of the diophantine equation \eqref{EQ1} is said to be primitive if $\gcd(A,B,C,D)=1.$
 \end{definition}
 Using  SAGE software \cite{sage}, we see that $Rank(E_6(\Bbb{Q}))=1$ and  $P=(-3,9)$ is  the generator of non-torsion subgroup of $E_6(\Bbb{Q})$.
 So, without taking into consideration of the inverse points, we know that every  point  of the form  $(X,Y)=n(-3,9)$ for some $n\in{\Bbb N}$, is also a non-torsion  point in  $E(\Bbb{Q})$.
 we have   $n(-3,9)=\left(\frac{\phi_n(-3,9)}{{\psi^2}_n(-3,9)} , \frac{\omega_n(-3,9)}{{{\psi^3}_n}(-3,9)}\right)$ where
 $\psi_n$ is the n-th division polynomial of $E_6$, $\phi_n=x{\psi^2}_n-{\psi_{n+1}}{\psi_{n-1}}$ and
 $\omega_n=(4y)^{-1}({\psi}_{n+2}{\psi^2}_{n-1}-{\psi}_{n-2}{\psi^2}_{n+1})$ (For more details see \cite{Wash} pages 81-84). Therefor, we can set $e={\psi}_n(-3,9)$, $b=\phi_n(-3,9)$ and $c=\omega_n(-3,9)$ in Eq \eqref{EQ4} to obtain a sequence of solution $(X_n,Y_n,U_n,V_n)$. For simplicity we omit $(-3,9)$ to obtain
\begin{equation}
\begin{array}{l}
  X_n={\phi_n}^4 + 1296{{\psi}_n}^8 + 864{\phi_n}{{\psi}_n}^6 + 72{\phi_n}^2{{\psi}_n}^4 + 144{\omega_n}{{\psi}_n}^5 \\
  \qquad\qquad - 24{\phi_n}^3{{\psi}_n}^2 + 4{\phi_n}^2{\omega_n}{{\psi}_n},\\
  \\
  Y_n=-864{\phi_n}{{\psi}_n}^6 - {\phi_n}^4 - 1296{{\psi}_n}^8 - 72{\phi_n}^2{{\psi}_n}^4 + 144{\omega_n}{{\psi}_n}^5\\
  \qquad\qquad + 24{\phi_n}^3{{\psi}_n}^2 + 4{\phi_n}^2{\omega_n}{{\psi}_n},\\
  \\
  U_n=4({\phi_n}^2 + 36{{\psi}_n}^4){\omega_n}{{\psi}_n},\\
  \\
  V_n=({\phi_n}^2 - 36{{\psi}_n}^4 - 12{\phi_n}{{\psi}_n}^2)({\phi_n}^2 - 36{{\psi}_n}^4 + 12{\phi_n}{{\psi}_n}^2).\\
\end{array}
\end{equation}
  Not all the solutions of $(X_n, Y_n, U_n, V_n)$ are primitive and some of them are multiples of $(21,19,20,7)$. For instance $(-3,9)$ leads to $(189,171,180,-63)=9(21,19,20,-7)$, which is not a new one.\\
  Let $(A_n,B_n,C_n,D_n)=(X_n/d_n, Y_n/d_n,U_n/d_n,V_n/d_n)$, in which $d_n=\gcd(X_n,Y_n,U_n,V_n)$. Obviously, $\{(A_n,B_n,C_n,D_n)\}$
  is a sequence of primitive solutions of diophantine equation \eqref{EQ2}. Using  SAGE, we computed
  $(A_n,B_n,C_n,D_n)$ for $2\leq n \leq 1000$ and presented some of new primitive solutions in the following.

\begin{align*}
 &n=2,\\
 &A_2=988521,&\\
 &B_2=-1661081,&\\
 &C_2=-336280,&\\
 &D_2=-1437599.\\
 \\
 &n=3,\\
 &A_3=-22394369951939,&\\
 &B_3=-59719152671941,&\\
 &C_3=-41056761311940,&\\
 &D_3=43690772126393.\\
 \\
 &n=4,\\
 &A_4=5009010521962601088594641,&\\
 &B_4=-959074737626305392403761,&\\
 &C_4=2024967892168147848095440,&\\
 &D_4=4156118808548967941769601.\\
 \\
 &n=5,\\
 &A_5=385103462588108468740542460457075040101,&\\
 &B_5=-58316597151277440454625613485820959901,&\\
 &C_5=163393432718415514142958423485627040100,&\\
 &D_5=-318497209829094206727124168815460900807.
  \end{align*}

\end{document}